\documentclass[12pt]{article}

\usepackage{newtxtext,newtxmath}
\usepackage{graphicx}

\usepackage[letterpaper,margin=1in]{geometry}

\usepackage{float}

\usepackage{ragged2e}
\usepackage{blindtext}

\renewenvironment{abstract}
	{\quotation}
	{\endquotation}

\date{}

\makeatletter
\renewcommand{\fnum@figure}{\textbf{Fig. \thefigure}}
\renewcommand{\fnum@table}{\textbf{Table \thetable}}
\makeatother

\usepackage{scicite}

\usepackage{url}

\def\scititle{
	The nexus between disease surveillance, adaptive human behavior and epidemic containment
}
\title{\bfseries \boldmath \scititle}

\author{
	Baltazar Espinoza$^{1\ast\dagger}$,
	Roger Sanchez$^{2\dagger}$,
	Jimmy Calvo-Monge$^{3\dagger}$
    Fabio Sanchez$^{2\dagger}$\and
	\small$^{1}$Biocomplexity Institute, University of Virginia, Virginia, USA.\and
	\small$^{2}$Escuela de Matemática-CIMPA, Universidad de Costa Rica, Ciudad Universitaria Rodrigo Facio,\and
    \small San José, 11501, Costa Rica.\and
    \small$^{3}$Centro de Investigación en Matemática Pura y Aplicada, Universidad de Costa Rica,\and
    \small Ciudad Universitaria Rodrigo Facio, San José, 11501, Costa Rica.\and
	\small$^\ast$Baltazar Espinoza. Email: baltazar.espinoza@virginia.edu\and
	\small$^\dagger$These authors contributed equally to this work.
}
\begin{document} 

% Insert the title and author list
\maketitle

\begin{abstract} \bfseries \boldmath
% 250 words max
Epidemics exhibit interconnected processes that operate at multiple time and organizational scales—a hallmark of complex adaptive systems. Modern epidemiological modeling frameworks incorporate feedback between individual-level behavioral choices and centralized interventions. Nonetheless, the realistic operational course for disease detection, planning, and response is often overlooked.
Disease detection is a dynamic challenge, shaped by the interplay between surveillance efforts and transmission characteristics. It serves as a tipping point that triggers emergency declarations, information dissemination, adaptive behavioral responses, and the deployment of public health interventions.
Evaluating the impact of disease surveillance systems as triggers for adaptive behavior and public health interventions is key to designing effective control policies.

We examine the multiple behavioral and epidemiological dynamics generated by the feedback between disease surveillance and the intertwined dynamics of information and disease propagation. 
Specifically, we study the intertwined dynamics between: $(i)$ disease surveillance triggering health emergency declarations, $(ii)$ risk information dissemination producing decentralized behavioral responses, and $(iii)$ centralized interventions.
Our results show that robust surveillance systems that quickly detect a disease outbreak can trigger an early response from the population, leading to large epidemic sizes. The key result is that the response scenarios that minimize the final epidemic size are determined by the trade-off between the risk information dissemination and disease transmission, with the triggering effect of surveillance mediating this trade-off.
Finally, our results confirm that behavioral adaptation can create a hysteresis-like effect on the final epidemic size.

% Each individual is responsible for their own decisions during an epidemic outbreak, people naturally exhibit a conservative, moderate, or risky attitude towards risk. This behavior is observed not only in finance, where it is more common to see and categorize these behaviors, but also from the onset of an epidemic up to a certain point in time. 
% In this manuscript, we study and contrast when behavioral change is enforced through a centralized policy, where the imposed rules must be strictly followed, versus emergency declarations made by public health officials that immediately affect the behavioral actions of each individual, that is, a decentralized policy. Through a behavioral awareness model, specific parameters are defined to represent the impact of centralization and decentralization, allowing an appropriate comparison. The trigger for the beginning of behavioral change is considered when the prevalence level reaches a pre-established threshold, $P^*$, set by health authorities.
\end{abstract}

\noindent

%Behavioral-epidemiology
Epidemics exhibit interconnected processes that operate at multiple time and organizational scales, a hallmark of complex adaptive systems~\cite{levin2003complex,levin1992problem,liu2007complexity}. The interdependency between human behavior and epidemics, influenced by endogenous feedback processes, is widely recognized. Together, behavioral-epidemiological dynamics generate a set of complexities shaping epidemic dynamics~\cite{bavel2020using,seale2020covid,moya2020dynamics,pagliaro2021trust,petherick2021worldwide,chen2018feedback,van2020using}.
At the individual level, transmission events are shaped by behaviors such as social interactions, mobility patterns, and interventions compliance~\cite{chinazzi2020effect,espinoza2020mobility,chang2021mobility}. At the social level, factors like healthcare infrastructure, government policies, and socioeconomic conditions play crucial roles in shaping epidemic progression and management~\cite{klepac2011synthesizing,glaubitz2024social,saad2023dynamics}.
Mathematical modeling has been effective in studying epidemic containment. There is an extensive literature of disease models incorporating distinct behavioral frameworks: threshold models of decision-making~\cite{qiu2022understanding,granovetter1978threshold,guilbeault2018complex,dodds2005generalized,morsky2023impact}, game theory models of behavioral choices~\cite{reluga2010game,bauch2004vaccination}, or behavioral choices based on cost-benefits trade-offs~\cite{espinoza2021asymptomatic,espinoza2022heterogeneous,espinoza2025impact,perrings2014merging,fenichel2011adaptive}.
However, the concurrent impacts of massive behavioral changes and the effects of interventions remain outstanding questions.

% Awareness
People do not passively comply with mandated interventions or suddenly change their behaviors. Information availability about the presence of a biological threat plays a key role in shaping the population's responses~\cite{funk2009spread,towers2015mass,schaller2011behavioural}.
Awareness of infection risk and associated illness costs create cost-benefit trade-offs that drive behavioral adoption~\cite{glaubitz2024social,chen2019imperfect}. 
Moreover, declarations made by public health officials lead to alterations in people's daily behavior over time, as knowledge and understanding of the infection risk change.
Consequently, preventive behaviors are typically adopted for a limited period, with potential relapses driven by pandemic fatigue and disease resurgence~\cite{petherick2021worldwide,nytfatigue,nyteurope}.
Information dissemination and disease contagion are coupled antagonistic spreading processes as they exert opposing effects on individual behaviors~\cite{de2016physics}. Each process exhibits its own dynamics, however, information spreading is inherently tied to the progression and detection of disease transmission.

% Biosurveillance
Conventional models assume the population's innate ability to adopt preventive behaviors and implement interventions arbitrarily during an epidemic. The realistic operational course that accounts for disease detection, planning, and response is often overlooked.
Disease surveillance systems are crucial for early warning, detection, and situational assessment of ongoing outbreaks. They provide valuable time to enhance our ability to respond to emerging biological threats~\cite{espinoza2023coupled,chen2024wastewater,margevicius2016biosurveillance}.
Disease detection is a dynamic challenge, shaped by the interplay between surveillance efforts and transmission characteristics. It serves as a tipping point that triggers emergency declarations, information dissemination, adaptive behavioral responses, and the deployment of public health interventions~\cite{perry2007planning,wolfe2021systematic}.
Incorporating the role of surveillance systems into mathematical modeling is critical to understanding the interdependence between epidemiological surveillance, disease dynamics and intervention strategies in an operational context.

\begin{figure}[H]
\centering
\includegraphics[width=\textwidth]{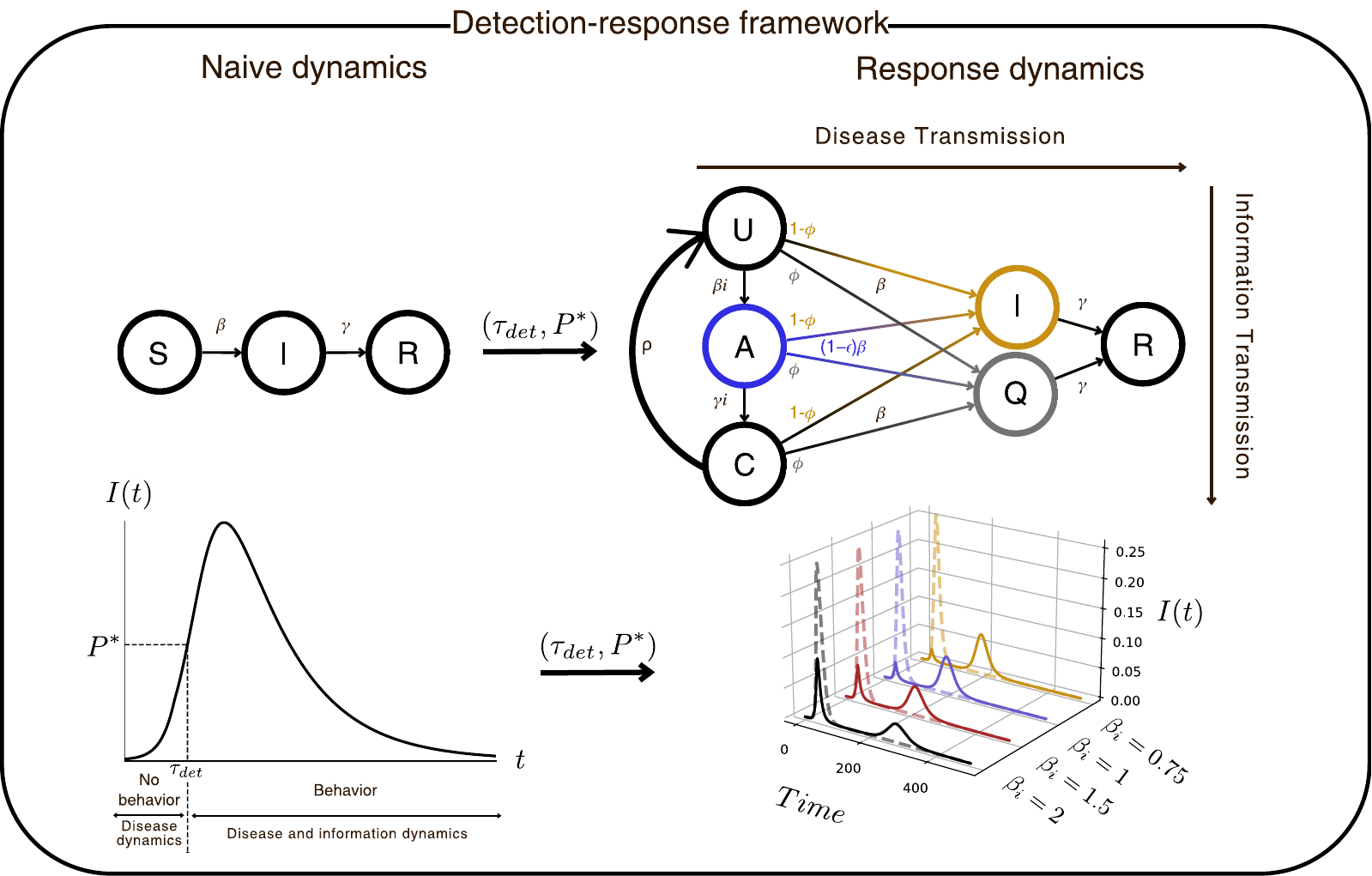}
\caption{
{\bf Schematic of our detection and response framework}
Our framework captures the intertwined dynamics between disease dynamics spreading over a naive population, defining the detection time, and the disease dynamics produced by including behavior and interventions. 
The detection and response times depend on the surveillance effort and disease progression. This alters the contagion dynamics by triggering dueling dynamics between information dissemination and disease transmission, which ultimately shape the effects of behavior and interventions on the epidemiological outcomes.
%
%Naive model with parameters $\beta$ and $\gamma$ which correspond to the infection rate to which susceptibles are exposed and to the rate at which infected individuals recover respectively. Moving on to the behavioral model, once the minimum prevalence level determined by $P$ has been reached, this adds three parameters that contribute to the transmission of information. The rates $\beta_i$, $\gamma_i$ and $\rho$ represent the horizontal transmission from the unaware to the aware as the infection rate, the recovery rate and which corresponds to the relapse in information transmission, respectively. The rate $\epsilon$ corresponds to the behavioral response strength and $\phi$ corresponds to the forced behavioral response.
% Illustrations show the time series of the infected individuals, respectively when the model is considered single shot and when is multiple shot. 
%Under a single shot framework, time series are illustrated to show how the disease progresses over time for different rates of information spread. The solid lines correspond to the behavioral model, while the dashed lines represent the naive model. Specifically, a maximum behavioral response strength is used, and infection parameters such that $R_0^d = 2.67$.
}
\label{fig:General_Model}
\end{figure}

We examine the multiple behavioral and epidemiological dynamics generated by the feedback between disease surveillance and the intertwined dynamics of information and disease propagation.
In this work, we study the intertwined dynamics between: $(i)$ disease surveillance triggering health emergency declarations, $(ii)$ risk information dissemination producing decentralized behavioral responses, and $(iii)$ centralized interventions.
Our modeling framework focuses on a realistic course of action based on testing, analysis, and responses. We assume the epidemic originates within a naive population unaware of the ongoing disease outbreak. A health emergency triggers awareness dissemination upon disease detection, temporarily changing individuals' understanding of infection risk. This alters their behavior and turns them into risk communicators.
We combine a probabilistic disease surveillance model with a behavioral-epidemiological model to analyze how the timing of health emergency declaration modulates behavioral responses and interventions, ultimately shaping epidemic outcomes.
We let health emergency declarations alter disease transmission by triggering information dissemination, which modifies behavioral dynamics, and the timing of interventions. Our modeling framework is schematically represented in Fig.~\ref{fig:General_Model}.
The SI Appendix provides a detailed formulation of the proposed modeling framework and a model extension incorporating awareness relapse.
We focus on the final epidemic size as a key epidemiological metric for the impact of the disease and information transmission.

Our results shed light on the intrinsic relationship between the surveillance effort, the disease basic reproductive number, and the risk information transmission inducing behavioral responses.
Specifically, we show that robust surveillance systems that quickly detect a disease outbreak can trigger an early response from the population, leading to high epidemic sizes under single-period responses (limited-duration behavioral adoption applied only once). In particular, the single-period most aggressive intervention is not always the best, commonly used strategy.
The key result is that the response scenarios that minimize the final epidemic size are determined by the trade-off between the risk information dissemination and disease transmission, with the triggering effect of surveillance mediating this trade-off.
On the other hand, our results align with previous findings that behavioral adaptation can create a hysteresis-like effect on the final epidemic size, a phenomenon not observed when only considering centralized interventions.

% The prevalence level of the disease is the point at which both models converge to enhance and strengthen the state of alert previously declared by health authorities. Therefore, the prevalence level of the disease is the primary and sole trigger for the state of awareness and disaggregation of susceptible individuals. Denoted by $P^*$, the triggering threshold for the activation of public awareness is identified by public health officials; this process consists of a awareness campaign transmitted through mass media, social networks, word of mouth, among others. \cite{azizi2020epidemics} Once this threshold is declared, it translates into terms of the number of quarantined and infected individuals at time $t^*$. When the naive model reaches a fixed prevalence value $I(t)=P^*$, the naive model becomes into behavioral.

% The risk associated with the disease is represented through various rates. This informational risk is transmitted to the population via different mechanisms triggered by the declaration of a health emergency by authorities. Declarations made by government officials lead to alterations in people's daily behavior over a period of time, as knowledge and understanding of the infection risk undergo changes. Experiencing a change in awareness and communicating to others about the existing infection risk is modeled through a rate that takes into account the economic and educational levels of individuals in the sample~\cite{azizi2020epidemics}.

\section*{Disease surveillance triggers behavioral responses and interventions deployment}

%The starting time of the population's behavioral response is given by the recognition of the epidemic, at the time the contagion reach a prevalence level of $P^*$ percent of the population, that is $t^*$ is such that $I(t^*)=P^*N$.

%NPIs can be absolutely necessary, particularly at the beginning of a pandemic due to the lack of knowledge and studies. For example, during the COVID-19 pandemic, some countries had very limited access to vaccines and in other cases, no access at all~\cite{wagner2021vaccine}.

\begin{figure}[H]
\centering
\includegraphics[width=\textwidth]{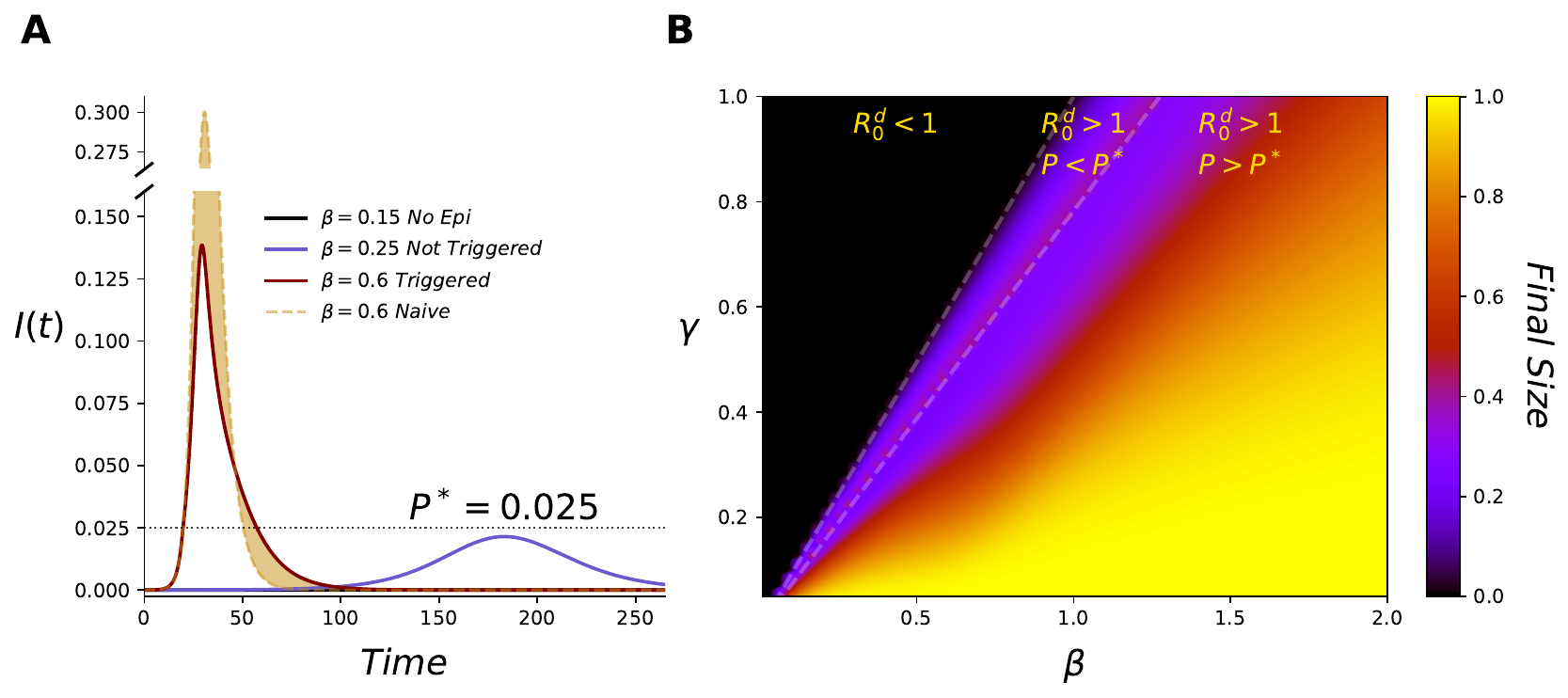}
\caption{{\bf The impact of disease transmission characteristics on the behavior-disease dynamics and the final epidemic size.}
{\bf (A)} Time series of infectious individuals showing the potential dynamics produced by the proposed framework. We assume an average infectious period $1/\gamma=10$ days, and $\mathcal{R}_0=\{0.75,1.25,3\}$, respectively.
{\bf(B)} Final epidemic sizes are attained by varying the disease infectiousness and period. We assume that the information dissemination exhibits a transmission likelihood of $\beta_i=1.5$, and that aware individuals show an average response period of $1/\gamma_i=10$ days, during which their susceptibility reduces by a factor $\epsilon=0.8$. We let the centralized response to quarantine a fraction $\phi=0.2$ of newly infected individuals and that surveillance effort triggers behavioral responses at a prevalence level of $P^*=0.025$.} 
\label{fig:timeSeries_R0}
\end{figure}

We begin by studying the impact of varying the disease basic reproductive number ($\mathcal{R}^d_0$) on triggering information and behavioral dynamics and the final epidemic size.
Fig.~\ref{fig:timeSeries_R0}(A-C) illustrate the diverse disease dynamics the proposed framework produces. Specifically, Fig.~\ref{fig:timeSeries_R0}A and B show the classical scenarios where the epidemic dies out if $\mathcal{R}^d_0<1$, and where the epidemic propagates when $\mathcal{R}^d_0>1$, without detection and in the absence of behavioral responses.
In contrast, Fig.~\ref{fig:timeSeries_R0}C shows that highly infectious diseases trigger behavioral responses that significantly reduce the peak size and prolong the epidemic for a given disease prevalence threshold.
Fig.~\ref{fig:timeSeries_R0}D illustrates the impact of varying the disease transmission rate and infectious period on the final epidemic size. Our results identify three regions (denoted by dashed lines) where the epidemic progression exhibits distinct dynamics: $(i)$ no epidemic propagation, where the epidemic dies out after onset ($\mathcal{R}^d_0<1$); $(ii)$ intermediate final epidemic sizes, where the epidemic propagates, producing low enough prevalence levels to avoid detection ($\mathcal{R}^d_0>1$); and $(iii)$ intermediate-high to high final epidemic sizes, where the epidemic reaches the detection threshold and triggers behavioral responses ($\mathcal{R}^d_0>>1$).

%The information contagion process makes individuals to be aware of the outbreak and take control measures on average during $10$ days, reducing their infection probability by $\epsilon=0.8$.

%Three time series are shown. When there is an epidemic and the behavioral phase is activated ((C), third area), parameters $\beta = 0.6$ and $\gamma = 0.2$ are used. When there is an epidemic, but the threshold of infected individuals is not surpassed, meaning there is no activation of behavior ((B), second area), where $\beta = 0.25$ and $\gamma = 0.2$. Finally, when the couple of parameters are over the blue zone, here there is no an epidemic, and thus there is no behavior ((A), third area) It is observed how $R_0^d$ decreases from left to right, showing decreasing severity as it progresses; and how the identity line, if we thought the picture as a Cartesian plane, divides epidemic and no epidemic.

\section*{Behavioral responses induce hysteresis in the final epidemic size}

In this section, we examine the impact of the competing social and biological dynamics on the final epidemic size. Specifically, we show that for single-period interventions: $(i)$ early responses do not minimize the final epidemic size, the impact of adaptive behavioral responses is influenced by the dynamics of risk information dissemination, where the final epidemic size exhibits a non-monotonic change as modifications in collective behavior accelerate; $(ii)$ Hysteresis in the final epidemic size is driven by decentralized behavioral responses. Centralized interventions alone show a monotonic final epidemic size as a function of disease detection ($ P*$), where behavioral response timing and stringency modulate the final size's non-monotonic change.

Hysteresis has previously been found in biological and social processes, such as drinking, where it was determined that the current state of drinkers depends on their previous state. In such a way that the susceptible-drinker and recovered-drinker flows contrast and help the drinking behavior developed~\cite{sanchez2007drinking}. %2025/03/05 Rg

% Different comparisons are made between the information and infection parameters, varying the basic reproductive number of the disease $R_0^d$ and the one of information $R_0^i$, in order to better visualize the relationship between both phases of coupling. In Fig. \ref{fig:timeSeries_R0} (D), three important areas can be observed. The first area corresponds to the darkest zone, which occurs when the infection parameters are such that $R_0^d < 1$, and therefore, an epidemic never occurs because the conditions are not sufficient for the disease to prevail. Right next to the first zone, there is the second area, mostly shown in gray and linearly delimited, the main characteristic of this area is that it satisfies $R_0^d > 1$; however, the model is never motivated to activate the behavioral phase of coupling, meaning $P < P^*$. Finally, there is the third area, which corresponds to the remaining area of the contour plot. The characteristics here are that $R_0^d > 1$ and additionally $P > P^*$, so an epidemic occurs and the behavioral phase of coupling is activated, reducing the epidemic final size compared to the prolonged first phase. 

\begin{figure}[H]
\centering
\includegraphics[width=\textwidth]{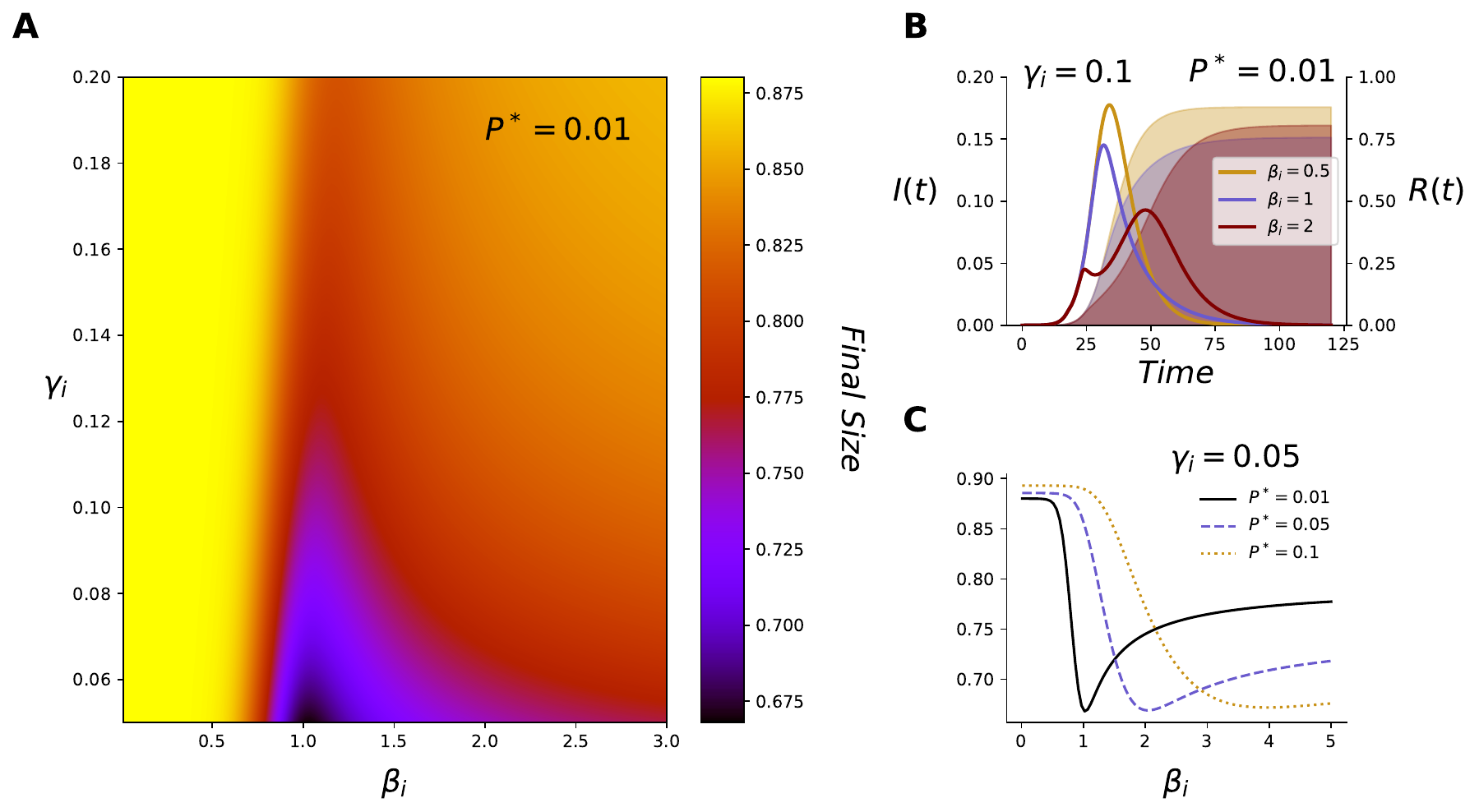}
\caption{
{\bf The impact of risk-information transmission characteristics on the final epidemic size.} 
{\bf (A)} Final size as a function of the risk-information transmission rate ($\beta_i$) and response withdrawal rate ($\gamma_i$). 
Our results show that, regardless of the behavioral response withdrawal rate, the final size exhibits a non-monotonic reduction as the risk-information transmission rate increases.
{\bf (B)} Time series of infected and recovered subpopulations for distinct risk-information transmission rates ($\beta_i = 0.5,1,2$), and response withdrawal rate $\gamma_i$ of $10$ days.
{\bf (C)} Final size as a function of the information transmission rate, for surveillance efforts leading disease detection at $P^*=0.01$, $P^*=0.05$, and $P^*=0.1$. 
% The earlier the disease transmission is detected, the lower the information transmission rate required to reduce the final epidemic size. Notice that, 
Early disease detection scenarios show a sharp minimum as a function of information transmission rate. Delayed disease detection requires higher information transmission rates and shows broader minimums.
%
%The parameter $\beta_i$ varies between $[0.01,3]$ while $\gamma_i$ varies between $[0.05,0.2]$. The remaining parameters, which are kept fixed during the simulation, are $\beta=0.6$, $\gamma=0.2$, $\epsilon=0.8$, $\phi=0.2$, $\rho=0$ (single shot) and $P=0.01$. {\bf (B)} A 3-dimensional perspective of (A) with the same parameters, for better understanding. {\bf (C)} Lowest minimum reached in simulation (A) ($\gamma_i=0.05$ with same set for $\beta_i$) along with other minima by varying the minimum prevalence level at which behavioral coupling is activated. Three simulations are done for three distinct prevalence thresholds $P=0.01,\ 0.05,\ 0.1$.
}
\label{fig:EFS_Information_R0}
\end{figure}

%\subsection*{Competing Social and Biological Dynamics Drive hysteresis}
% regardless of the behavioral response period ($1/\gamma_i$)
We find that changes in the information dissemination rate ($\beta_i$), mediate the timing of behavioral responses that produce a minimum on the final epidemic size. Our results in Fig.~\ref{fig:EFS_Information_R0}(A-B) show that behavioral responses triggered too soon produce large epidemics exhibiting reemergence dynamics. On the other hand, delayed behavioral responses are triggered at prevalence levels such that the final size is effectively reduced. However, if the behavioral responses are significantly delayed, the epidemic may have already reached a critical point to produce large outbreaks.
Moreover, Fig.~\ref{fig:EFS_Information_R0}C shows that disease surveillance efforts influence the information dissemination rate minimizing the final epidemic size ($P^*$).
Our results highlight the impact of the intertwined nature of the dueling social and biological dynamics. We show that behavioral responses non-monotonically impact the final epidemic size, depending on the adoption rate of collective behavior and the disease detection time. Furthermore, the rapid dissemination of risk information can alleviate the lack of robust surveillance systems. In such scenarios, delayed disease detection is compensated by the swift spread of awareness.

%The impact of the behavioral phase of the model on the final size of the epidemic can also be observed in Fig. (\ref{fig:EFS_Information_R0}). This figure represents how much the impact of the epidemic can be reduced if people adopted better behaviors or if they were more influenced by other individuals who have already adopted certain behaviors and are paying attention to the ongoing epidemic. The area of interest, of course, corresponds to where the lowest values are reached. In this case, it is seen that the lower values are obtained for a very particular region of $\beta_i$, which is less dependent on $\gamma_i$. This suggests that the ease of information transmission prevails over how long people remain in a state of alert. Of course, the lower the gamma, the smaller the final size of the epidemic, as the contagion is favored by the greater number of individuals carrying the information. An extension of this figure is shown in the appendix, where more scenarios are considered, providing a greater understanding of this effect.

The intuition behind the non-linear impact of information dissemination on the final epidemic size is represented in Fig.~\ref{fig:EFS_Information_R0}B. We show selected time series of the infected and recovered subpopulations, for $\beta_i$ values that produce distinct behavioral-disease dynamics. The risk-information dissemination rate producing the lowest number of recovered individuals corresponds to the intermediate value. Moreover, lower (greater) $\beta_i$ values generate a very delayed (a very early) behavioral response, producing a minimal impact on the epidemic size. 

%or so high that it initially reduces the epidemic but leaves a residual risk factor, increasing the final size in a second wave. This analysis implicitly interacts with the disease surveillance effort ($P^*$).

% \begin{figure}[H]
% \centering
% \includegraphics[width=\textwidth]{EFS_Information_R0_TS.pdf}
% \caption{{\bf Representation of risk-information transmission as time series.} Time series of infection and recovery compartments for three risk-information transmission rates ($\beta_i$ = 0.5,1,2) and response withdrawal rate $\gamma_i$ of $10$ days.}
% \label{fig:EFS_Information_R0_TS}
% \end{figure}

%The calibration of behavioral parameters is of great importance, as in most cases they are the result of disease surveillance, and thus they must reflect the consequences of raising alerts either early or late. There are people who comply with the interventions previously stipulated by health authorities; however, there are also those who go against the decisions made. In~\cite{chen2024simple}, an analysis is carried out in which vaccination and non-pharmaceutical interventions are considered, showing how the correlation between individuals' behaviors towards these two strategies affects epidemiological outcomes. Through simulations, the positive and negative consequences are shown based on the parameters; in particular, a reduction in the epidemic burden is observed due to interventions modulated according to disease surveillance.

 % Compare between a centralized vs decentralized policy % Present and highlight if people keep in aware compartmental

\begin{figure}[H]
\centering
\includegraphics[width=\textwidth]{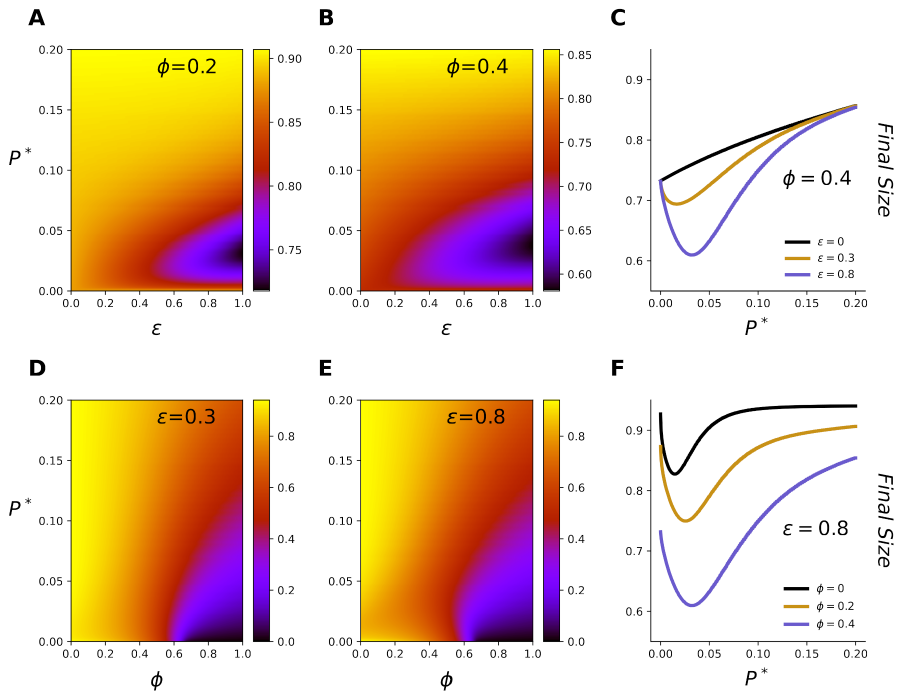}
\caption{
{\bf The trade-off between centralized interventions and adaptive behavior, on the final epidemic size.} 
{\bf (A and B)} Final epidemic sizes as a function of the behavioral response stringency ($\epsilon$) and the surveillance efforts ($P^*$), for quarantine capacities of $\phi=0.2$ and $\phi=0.4$. Our results show non-monotonic changes on the final epidemic size as the surveillance effort decreases. {\bf (C)} The local minimum emerges and is exacerbated by increments on the behavioral response stringency.
{\bf (D and E)} Final epidemic sizes as a function of quarantine capacity ($\phi$) and surveillance efforts ($P^*$), for behavioral response stringency of $\epsilon=0.3$ and $\epsilon=0.3$. Low stringency of behavioral responses ($\epsilon=0.3$) vanishes the non-monotonic changes on the final epidemic size. High stringency of behavioral responses ($\epsilon=0.8$) exhibit non-monotonic changes on the final epidemic size regardless of the quarantine capacity.
%
%Different scenarios by modifying the parameter $\phi=0.2,\ 0.4$ respectively for (A) and (B) in order to get the epidemic final size, where $\epsilon$ is varied between $[0,1]$ and $P$ between $[0,1]$, while keeping the rest of the parameters fixed, $\rho=0$ (single shot), $\beta =0.6$, $\gamma =0.2$, $\beta_i=1.5$ and $\gamma_i=0.1$. {\bf (D and E)} Images representing the same as in (A and B), with the difference being that the epidemic final size is obtained by varying $\phi$ instead of $\epsilon$, and considering two cases for $\epsilon=0.3,\ 0.8$ respectively for (D and E). {\bf (C and F)} Vertical trajectories obtained from (B and E) for (C and F) respectively.
}
\label{fig:BehaviorParameters_grids}
\end{figure}

We next explore the effects of decentralized adaptive behavioral responses of aware individuals ($\epsilon$) and centralized intervention quarantining a fraction of the newly infected individuals ($\phi$), for distinct disease surveillance efforts ($P^*$).
Fig.~\ref{fig:BehaviorParameters_grids}(A-C) show that while both centralized and decentralized responses decrease the final epidemic size, non-monotonic changes are exhibited only in the presence of adaptive behavioral responses. As behavioral efforts decrease ($\epsilon=0$), the final epidemic size shows monotonic increments. However, as behavioral responses become more prominent, the final epidemic size exhibits non-monotonic increments.
Furthermore, our results in Fig.~\ref{fig:BehaviorParameters_grids}(D-F) show that the non-monotonic effect on the final size exhibits limited sensitivity to changes in quarantine capacity. While increases in quarantine capacity amplify the impact of behavioral responses on reducing the disease burden, this does not affect the trade-off between surveillance and behavioral responses. In these scenarios, delayed disease detection combined with strong behavioral responses minimizes the final epidemic size.

%The impact on the final epidemic size of starting behavioral change at different times is explored in the following set of simulations. Fig. \ref{fig:BehaviorParameters_grids} (C) shows the impact of behavioral response strength ($\epsilon\in[0,1]$), by reducing the infectious rate for aware individuals during its period of awareness, on the final size for different prevalence thresholds. The minimum final size is reached almost at same prevalence levels ($P^*$) for different values of response strength ($\epsilon$). The higher the population response the lower the final size. The impact is maximized at a particular $P^*$. Moreover, Fig. \ref{fig:BehaviorParameters_grids} (F) shows the impact of individuals self-quarantine decision after being infected. Since infected quarantined individuals are assumed not infectious, removing a fraction $\phi$ of infected individuals from the ``effective'' infectious population has a dramatic impact on the reduction of the final epidemic size. The prevalence value at which the minimum final size is reached seems to be slightly sensitive to the self-quarantine levels. The sensitivity to changing the parameter $\phi$ is more notable between panels (A) and (B) compared to the change in $\epsilon$ in panels (D) and (E). However, panels (A,B) shows how the same pattern emerges when $P^*$ has a very low value, indicating that raising early alerts is counterproductive. 

\section*{Conclusion}

Disease surveillance systems offer early warnings of outbreaks before they escalate into public health emergencies.
The detection of biological threats represents trigger events that start a variety of signaling cascades that activate multiple responses.
The dissemination of an alert status during epidemics often exhibits transmission dynamics independent of the disease progression, and can show multiple transmission dynamics~\cite{towers2015mass,margevicius2016biosurveillance,funk2009spread,funk2010endemic}. Social media and news, or mouth-to-mouth information dissemination, could cause massive awareness, potentially leading to collective changes in individuals' behaviors.
In this work, we assume that the implementation of epidemic interventions and the adoption of risk-avoidance behaviors are preceded by identifying a biological threat. Our coupled disease surveillance model serves as an operational tipping point for public health authorities and the population's spread of awareness leading to preventive behaviors.

% Quick detection and slow information flow balance delayed detection and quick information flow.

%%%%%%%%%%%%Fabio%%%%%%%%%%%%%%%%%%
Our results highlight the role of disease surveillance as a trigger for public health interventions and adaptive human behavior, emphasizing the complex feedback mechanisms that shape the dueling dynamics of disease progression and risk awareness. We showed that surveillance efforts, by influencing the timing of emergency declarations and triggering risk communication and behavioral responses, significantly affect the final size of the epidemic.
Moreover, our findings underscore the importance of surveillance systems in modulating the effectiveness of interventions to mitigate the impact of epidemics. Rapid detection of biological threats enhances public awareness, prompting early behavioral changes and the implementation of intervention strategies. However, our results reveal that the timing of behavioral adoption is critical, as both early and delayed responses can lead to suboptimal outcomes. 
In particular, we identify a trade-off between risk information dissemination and disease transmission dynamics, where an overly aggressive early response may not always result in the smallest epidemic size. This result is crucial for strategic planning in determining the optimal timing for triggering behavioral responses and interventions.
Furthermore, our analysis illustrates how decentralized behavioral adaptations, influenced by social awareness and information spread, can create hysteresis-like effects on the final epidemic size. The emergence and modulation of this phenomenon in the presence of behavioral responses suggests that the effectiveness of public health policies heavily depends on the dynamic interplay between individual decision-making and collective responses. Our results align with existing literature showing that individual-specific control measures outperform population-wide measures~\cite{lloyd2005superspreading,qiu2022understanding,morsky2023impact,mahmud2024adaptive}.

Our extended model, incorporating awareness relapse, provides further insights into the dynamics of epidemic control. The potential for individuals to revert to pre-awareness behaviors suggests that sustained public health messaging and reinforcement strategies are necessary to prevent a resurgence of disease transmission. Additionally, our results align with previous studies indicating that effective disease containment requires balancing top-down intervention policies and bottom-up behavioral responses~\cite{espinoza2023coupled,dixit2023airborne,glaubitz2024social}.
Finally, this work provides a quantitative framework for understanding the interplay between disease surveillance, public health decision-making, and the joint dynamics of adaptive human behavior and epidemics. The insights gained from this study aim to inform the design of effective epidemic containment strategies by optimizing the timing of emergency declarations, integrating the impact of risk communication, and accounting for the dynamics of behavioral responses to minimize the societal burden of infectious diseases.

% 

%%%%%%%%%%%%%%%% REFERENCES %%%%%%%%%%%%%%%

\clearpage % Clear all remaining figures and tables then start a new page

% The list of references goes after the main text and before the acknowledgements
% When preparing an initial submission, we recommend you use BibTeX, like this:

%\bibliography{science_template} % for a file named science_template.bib
%\bibliographystyle{sciencemag}

% After the paper has completed peer review and been revised ready for acceptance,
% you should comment out the lines above and copy-paste the contents of your .bbl
% file here instead. This will help ensure that our conversion software works correctly.
% Remember to re-run BibTeX first - check the timestamp!
%
% Example of the first three entries copy-pasted from science_template.bbl:
%
%\begin{thebibliography}{1}
%
%\bibitem{example}
%A.~N. {Author}, An example reference. \emph{Journal of Improbable Research}
%  \textbf{1}, 67 (2020).
%
%\bibitem{example2}
%F.~M. {Surname}, S.~{Author}, A second example. \emph{Interesting Research
%  Letters} \textbf{32}, 897 (2019).
%
%\bibitem{example_preprint}
%P.~{One}, P.~{Two}, P.~{Three}, {An unpublished preprint}. \emph{preprint}
%  (2021), arXiv:2101.12345.
%
%\end{thebibliography}

%%%%%%%%%%%%%%%% ACKNOWLEDGEMENTS %%%%%%%%%%%%%%%

\section*{Acknowledgments}

%%%%%%%%%%%%%%%% SUPPLEMENT LIST %%%%%%%%%%%%%%%

% List the contents of your Supplementary Materials, including the numbers of any
% supplementary figures, tables, external data files etc. and any references that are
% cited only in the supplement. In this example, refs. 7-8 are cited only in the supplement.
% Fill out your numbers accordingly and delete any lines that aren't applicable.
\subsection*{Supplementary materials}
Materials and Methods\\
Supplementary Text\\
Figs. S1 to S9\\
References \textit{(7-\arabic{enumiv})}\\

%%%%%%%%%%%%%%%% END OF MAIN TEXT %%%%%%%%%%%%%%%

\newpage

%%%%%%%%%%%%%%%% START OF SUPPLEMENT %%%%%%%%%%%%%%%

% Figures, tables, equations and pages in the supplement are numbered S1, S2 etc.
\renewcommand{\thefigure}{S\arabic{figure}}
\renewcommand{\thetable}{S\arabic{table}}
\renewcommand{\theequation}{S\arabic{equation}}
\renewcommand{\thepage}{S\arabic{page}}
\setcounter{figure}{0}
\setcounter{table}{0}
\setcounter{equation}{0}
\setcounter{page}{1}

%%%%%%%%%%%%%%%% SUPPLEMENT TITLE PAGE %%%%%%%%%%%%%%%

\begin{center}
\section*{Supplementary Materials for\\ \scititle}

% Author list for the supplement
% Indicate the corresponding authors, but do NOT include institutions here
% It would be nice if the template auto-generated this, but doing so is complicated...
	Baltazar Espinoza$^{1\ast\dagger}$,
	Roger Sanchez$^{2\dagger}$,
	Jimmy Calvo-Monge$^{3\dagger}$
    Fabio Sanchez$^{2\dagger}$\\

	% Identify at least one corresponding author, with contact email address
	\small$^\ast$Baltazar Espinoza. Email: be8dq@virginia.edu\\
	% Joint contributions can be indicated like this
	\small$^\dagger$These authors contributed equally to this work.
\end{center}

% Fill out the numbers for each type of supplementary material,
% and delete any lines that aren't applicable.
% These are just example numbers that don't match the rest of this template.
\subsubsection*{This PDF file includes:}
Materials and Methods\\
Supplementary Text\\
Figures S1 to S9\\

% \subsubsection*{Other Supplementary Materials for this manuscript:}

\newpage

%%%%%%%%%%%%%%%% MATERIALS AND METHODS %%%%%%%%%%%%%%%

\subsection*{Materials and Methods}

\subsubsection*{Coupling disease surveillance and epidemic dynamics}

Our coupled behavioral-epidemiological model depicted in Figure~\ref{fig:General_Model} assumes the epidemic initially spreads among a naive population with a Susceptible-Infected-recovered scheme formalized by equation~\eqref{eq:naive}
\begin{equation}
\begin{split}
\dot{S} &= -\beta S \frac{I}{N},\\
\dot{I} &= \beta S \frac{I}{N}-\gamma I,\\
\dot{R} &= \gamma I,
\end{split}
\label{eq:naive}
\end{equation}
where $\beta$ corresponds to the transmission rate, and $\gamma$ is the recovery rate. 
We couple the disease progression model with a surveillance model formulated using a probabilistic approach to compute the expected time required to identify at least a single infected individual using power calculations.
Assuming a homogeneous mixing population and a random sampling strategy, we compute the cumulative probability of detecting at least one infectious individual.
The proposed disease detection model incorporates the dynamic progression of the epidemic into the daily probability of detecting infectious individuals.
The probability of detecting at least a single positive case by time $t$, for a disease prevalence level of $\phi_s$ on day $s$, for surveillance effort of $N_s$ daily tested individuals, is given by
\begin{equation}
\mathbb{P}_{det}= 1 - \prod_{s=1}^{t} (1-\phi_s)^{N_s}.
\end{equation}
Finally, we let the disease detection time ($\tau_{det}$) to be reached at the time at least a single infectious individual is detected with at least 95\% confidence, $\mathbb{P}_{det}=0.95$.
The detection time determines the disease prevalence threshold ($P^*$) at which health emergency declarations are issued. Note that $P^*$ is not arbitrarily selected; instead it depends on the preselected surveillance effort.

\subsubsection*{Disease surveillance triggers interventions, risk information spread and adaptive human behavior}

Upon health emergency declaration, we assume the onset of centralized intervention quarantining a fraction $\phi$ newly infected individuals ($Q$), and that the susceptible population is divided into the following subgroups based on their knowledge about the epidemic:
Unaware ($U$), who are naive of the infection risk; Aware ($A$), who are informed of the infection risk, adopt preventive behaviors, and are ready to convey information; and Careless ($C$), who are informed of the infection risk but do not adopt preventive behaviors or convey information.
Dissemination of awareness (information about infection risk) begins with a single aware individual conveying the information to unaware individuals at a transmission rate $\beta_i$. Aware individuals adopt preventive behaviors that reduce their susceptibility to disease transmission by a factor of $(1-\epsilon)$ during an average period of $1/\gamma$ days. After this period, individuals become careless and stop adopting preventive behaviors, exhibiting the same susceptibility to disease transmission as unaware individuals.
Moreover, we assume the health emergency declaration triggers quarantine as a mandated intervention. Regardless of their knowledge status, a fraction $\phi$ of the newly infected individuals are quarantined until recovery, during which they cannot transmit the disease.
The dynamics as mentioned above are formalized by the set of differential equations in Equation~\eqref{eq:Single_Behavior}
%For the second part of the coupling, two submodels are considered. The first model, also referred to as single shot in the main text, disaggregates the susceptible category into three categories: unaware, aware, and careless. Additionally, it considers infected individuals and includes the quarantined compartment, while the recovered category remains unchanged. The particularity of this submodel lies in the fact that there is no more than one loop in the transmission of information. The following system of equations describes the submodel:
%
\begin{equation}
\begin{split}
\dot{U} &= -\beta_i U \frac{A}{N}-\beta U \frac{I}{N},\\
\dot{A} &= \beta_i U \frac{A}{N}- (1-\epsilon)\beta A \frac{I}{N}-\gamma_i A,\\
\dot{C} &= \gamma_i A-\beta C \frac{I}{N},\\
\dot{Q} &= \phi\beta\left[ U+C+(1-\epsilon) A \right]\frac{I}{N}-\gamma Q,\\
\dot{I} &=  (1-\phi)\beta\left[ U+C+(1-\epsilon) A \right]\frac{I}{N}-\gamma I,\\
\dot{R} &= \gamma (Q + I),
\end{split}
\label{eq:Single_Behavior}
\end{equation}
where $\epsilon$ represents the stringency of the behavioral response, $\phi$ denotes the proportion of quarantined individuals, $\beta_i$ is the information transmission rate, $\gamma_i$ correspond to the information transmission and recovery rates, and $\beta$ and $\gamma$ are the disease transmission and recovery rate, respectively.

\subsubsection*{Adaptive behavior relapse model}

We extend the model~\eqref{eq:Single_Behavior} by incorporating a potential relapse of careless individuals that allows for multiple adoption cycles of precautionary behavior. We assume that a careless individual ($C$) loses interest in the epidemic progression and returns to unaware status ($U$) at a rate $\rho$, after which it could become aware again by interacting with aware individuals.
The dynamics mentioned above are represented in model~\eqref{eq:Multiple_Behavior}
\begin{equation}
\begin{split}
\dot{U} &= -\beta_i U \frac{A}{N}-\beta U \frac{I}{N} + \rho C,\\
\dot{A} &= \beta_i U \frac{A}{N}- (1-\epsilon)\beta A \frac{I}{N}-\gamma_i A,\\
\dot{C} &= \gamma_i A-\beta C \frac{I}{N} -\rho C,\\
\dot{Q} &= \phi\beta\left[ U+C+(1-\epsilon) A \right]\frac{I}{N}-\gamma Q,\\
\dot{I} &=  (1-\phi)\beta\left[ U+C+(1-\epsilon) A \right]\frac{I}{N}-\gamma I, \\
\dot{R} &= \gamma (Q + I).
\end{split}
\label{eq:Multiple_Behavior}
\end{equation}

\subsection*{Scenario robustness}

In this section, we explore the robustness of our results as a function of distinct disease and information scenarios. Specifically, we study the impact of varying the length of the infectious and behavior periods and disease and information transmission likelihood on modulating the final epidemic size.

\subsubsection*{Single period behavioral response scenario}

% Final epidemic size as a function of disease transmission dynamics, for distinct information dynamic scenarios

% \begin{figure}[H]
% \centering
% \includegraphics[width=\textwidth]{Single_Shot.pdf}
% \caption{{\bf Single shot model flowchart.} Simplified model when it is such that the dynamics do not consider a relapse, meaning there is no return from the careless to the unaware state. Compared to the multiple shot scenario in the main text $\rho=0$.}
% \label{fig:Single_Shot}
% \end{figure}

\begin{figure}[H]
\centering
\includegraphics[width=\textwidth]{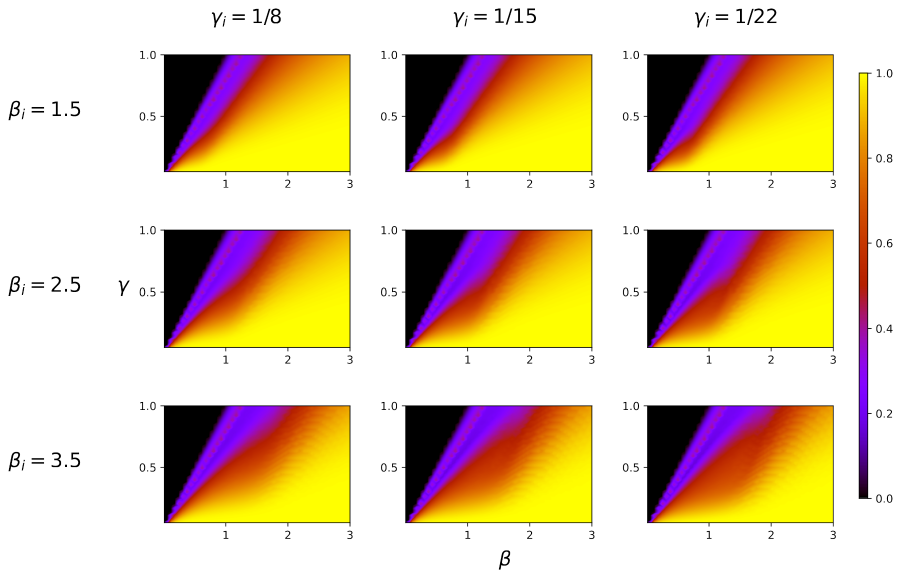}
\caption{{\bf Grid of informative parameters comparison.} Shows the gradual change when different infection and recovery rates are set in the informative section of the model, providing an idea of the 4-dimensional motion presented through the model's sensitivity. The grid is composed by the sets of $\beta_i\in\{1.5,2.5,3.5\}$ and $\gamma_i\in\{1/8,1/15,1/22\}$; while the other parameters are considered fixed.}
\label{fig:Grid_InformativeParametersComparison}
\end{figure}

\begin{figure}[H]
\centering
\includegraphics[width=\textwidth]{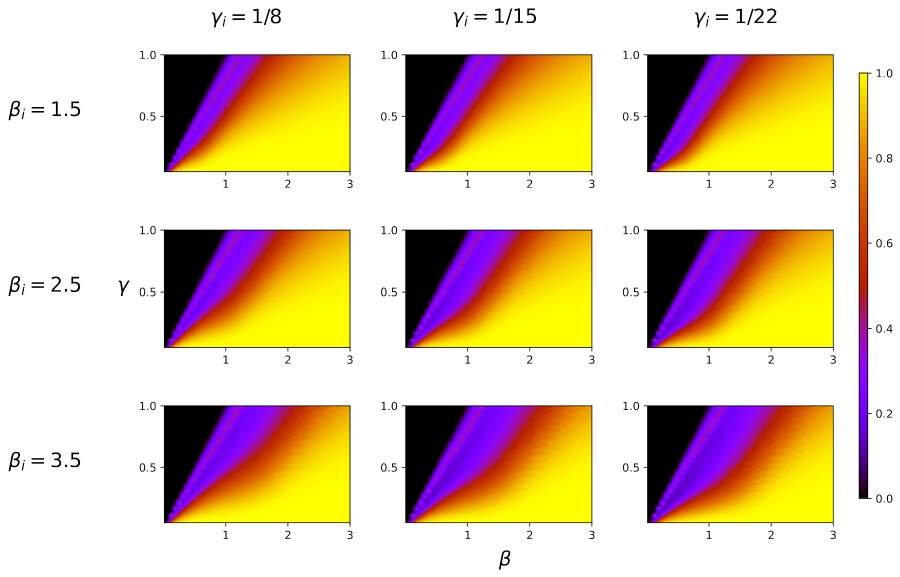}
\caption{{\bf Grid of informative parameters comparison in a multiple shot framework.} Shows the gradual change when different infection and recovery rates are set in the informative section of the model, providing an idea of the 4-dimensional motion presented through the model's sensitivity. The grid is composed by the sets of $\beta_i\in\{1.5,2.5,3.5\}$ and $\gamma_i\in\{1/8,1/15,1/22\}$; while the other parameters are considered fixed with $\rho=1/20$.}
\label{fig:Grid_InformativeParametersComparison_MultipleShot}
\end{figure}

\begin{figure}[H]
\centering
\includegraphics[width=\textwidth]{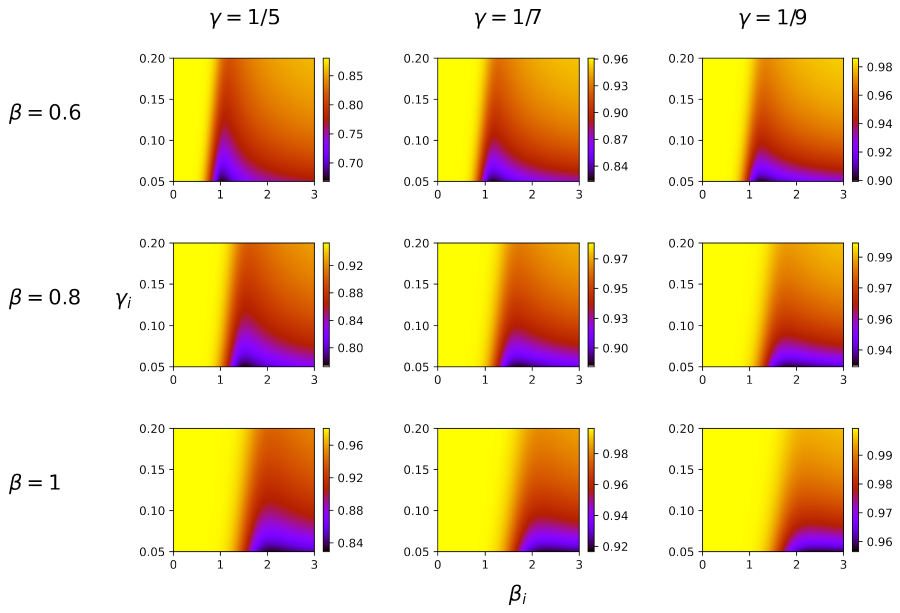}
\caption{{\bf Grid of infectious parameters comparison.} Shows the gradual change when different infectious and recovery rates are set in the infectious section of the model. Each subplot has its color bar showing the distribution of values. The grid is composed by the sets of $\beta\in\{0.6,0.8,1\}$ and $\gamma\in\{1/5,1/7,1/9\}$; while the other parameters are considered fixed.}
\label{fig:GridInfectionParametersComparison}
\end{figure}

\begin{figure}[H]
\centering
\includegraphics[width=\textwidth]{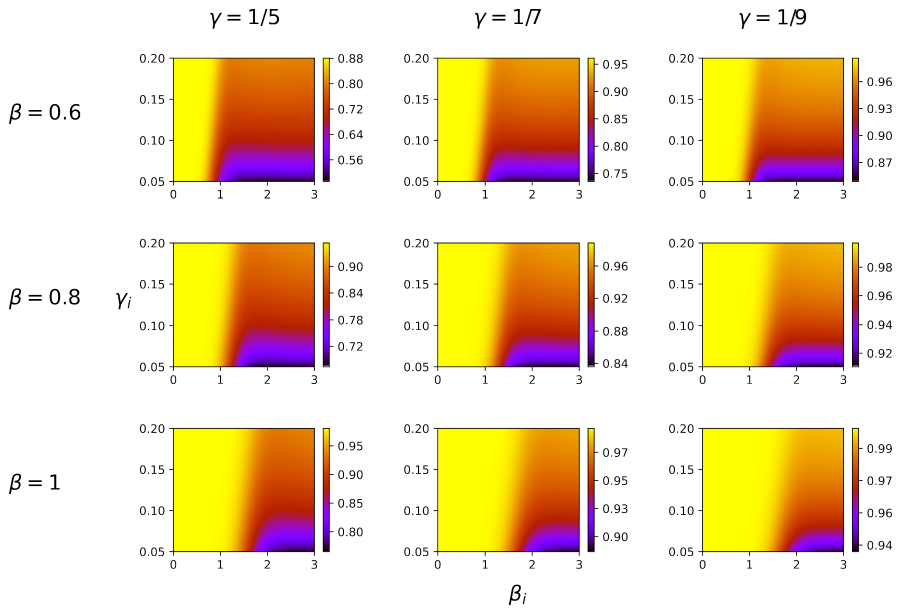}
\caption{{\bf Grid of infectious parameters comparison in a multiple shot framework.} Shows the gradual change when different infectious and recovery rates are set in the infectious section of the model. Each subplot has its own color bar showing the distribution of values. The grid is composed by the sets of $\beta\in\{0.6,0.8,1\}$ and $\gamma\in\{1/5,1/7,1/9\}$; while the other parameters are considered fixed with $\rho=1/20$.}
\label{fig:GridInfectionParametersComparison_MultipleShot}
\end{figure}

\begin{figure}[H]
\centering
\includegraphics[width=\textwidth]{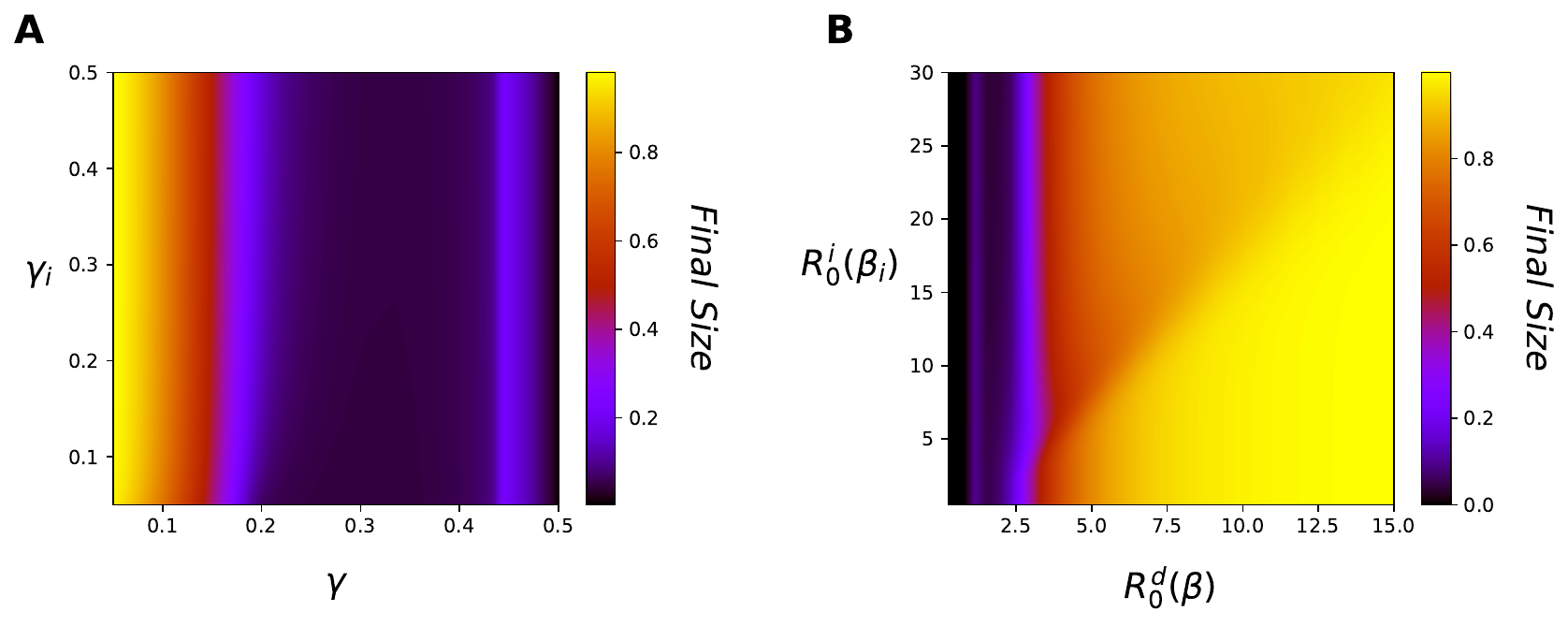}
\caption{{\bf Comparison between rates of the same category in a single shot.} Contrast between behavioral and informational parameters represented with a heatmap, accompanied by a 3-dimensional perspective for better visualization. Recovery rates vary between $[0.05,0.5]$ while contagion rates vary between $[0.05,3]$. The comparison of recovery rates is plotted in {\bf (A)}. Similarly, the contagion rates are contrasted in {\bf (B)}.}
\label{fig:ComparisonSameType}
\end{figure}

\begin{figure}[H]
\centering
\includegraphics[width=\textwidth]{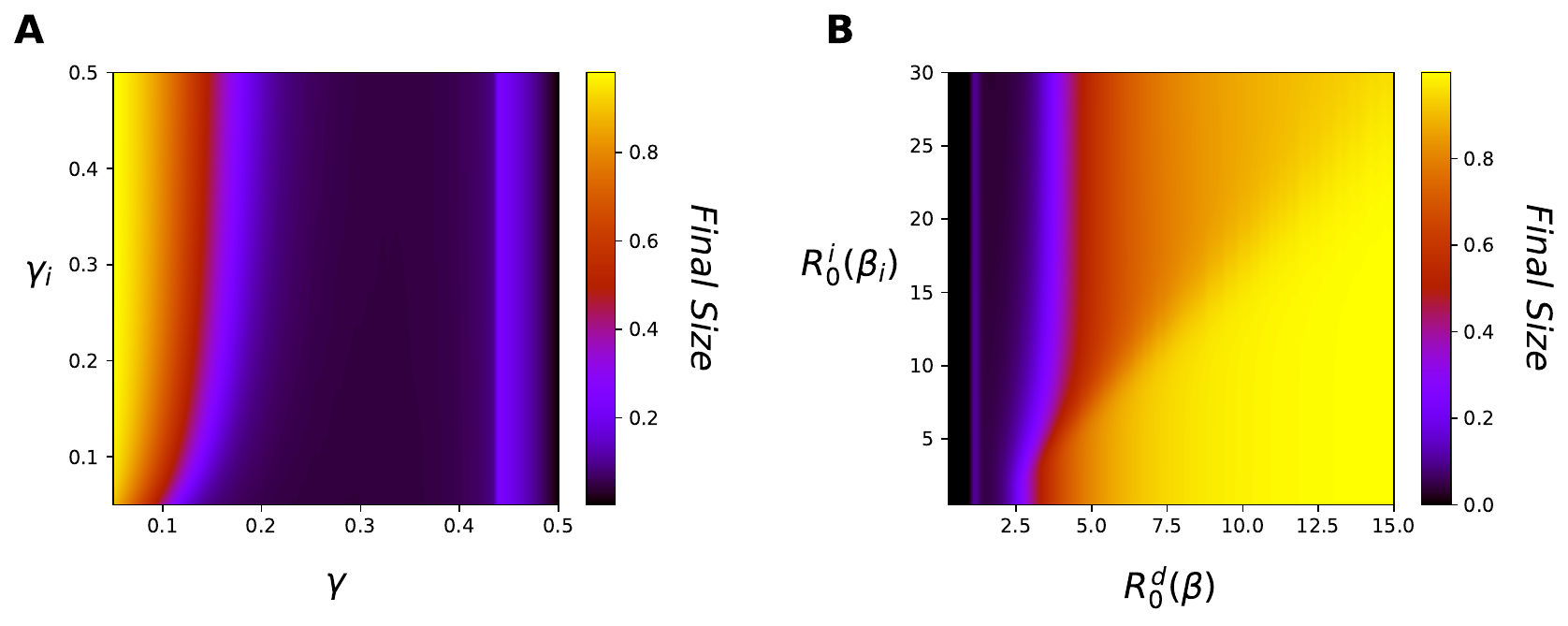}
\caption{{\bf Comparison between rates of the same category in a multiple shot.} Contrast between behavioral and informational parameters represented with a heatmap, accompanied by a 3-dimensional perspective for better visualization. Recovery rates vary between $[0.05,0.5]$ while contagion rates vary between $[0.05,3]$. The comparison of recovery rates is plotted in {\bf (A)}. Similarly, the contagion rates are contrasted in {\bf (B)}. }
\label{fig:ComparisonSameType_MultipleShot}
\end{figure}

\subsubsection*{Behavioral response with relapse scenario}

% \begin{figure}[H]
% \centering
% \includegraphics[scale=0.75]{Disease_Contagion_MS.pdf}
% \caption{{\bf Several scenarios according to the severity of information spread.} Under a single shot framework, time series are illustrated to show how the disease progresses over time for different rates of information spread. The solid lines correspond to the behavioral model, while the dashed lines represent the naive model. Specifically, a maximum behavioral response strength is used, and infection parameters such that $R_0^d = 2.67$.}
% \label{fig:Disease_Contagion_MS}
% \end{figure}

\begin{figure}[H]
\centering
\includegraphics[width=\textwidth]{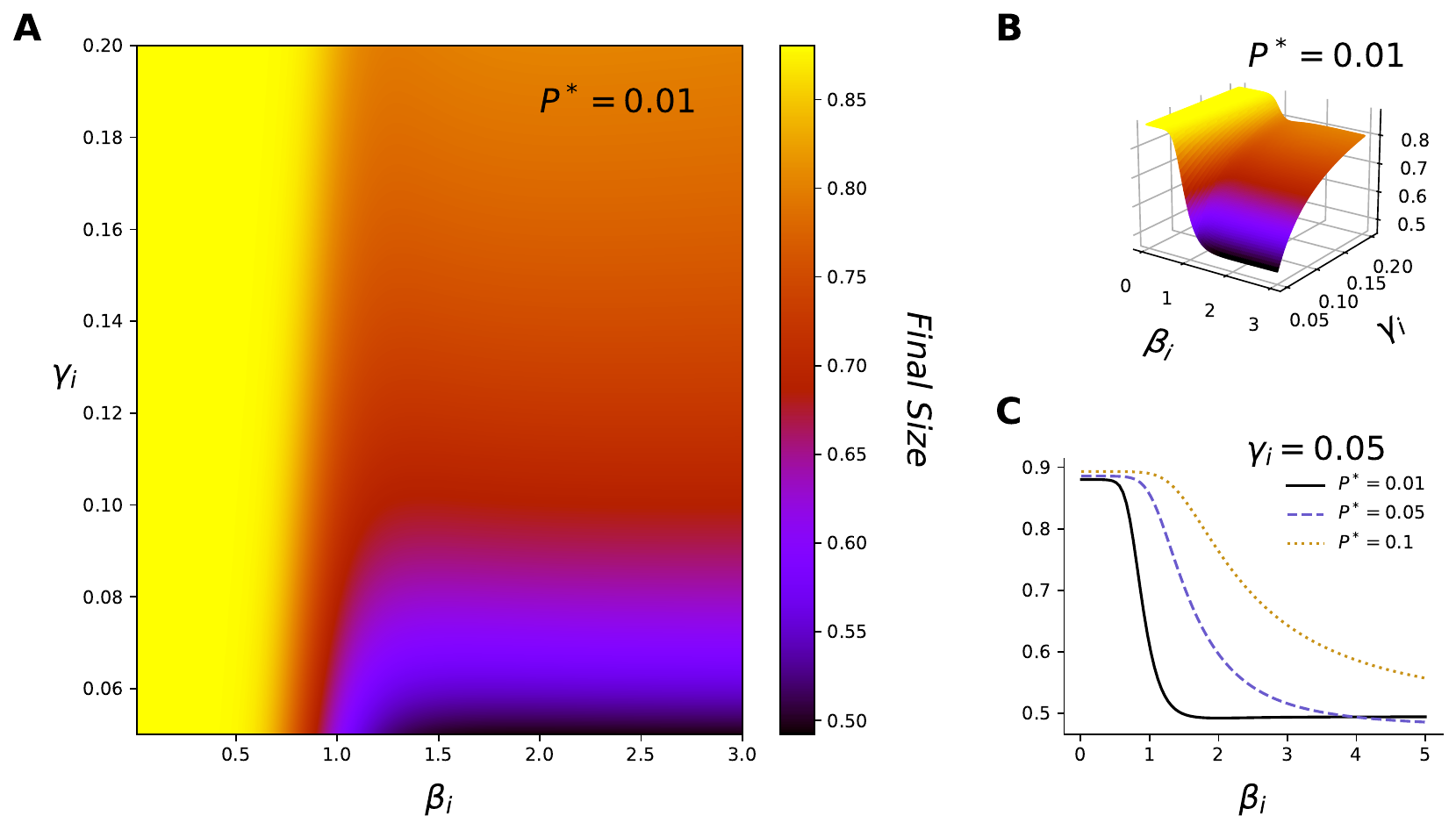}
\caption{{\bf Impact of behavior in the epidemic final size based on comparison of informative parameters in a multiple shot framework.} {\bf (A)} Heatmap of the epidemic final size based on the informative parameters. The parameter $\beta_i$ varies between $[0.01,3]$ while $\gamma_i$ varies between $[0.05,0.2]$. The remaining parameters, which are kept fixed during the simulation, are $\beta=0.6$, $\gamma=0.2$, $\epsilon=0.8$, $\phi=0.2$, $\rho=1/20$ (multiple shot) and $P=0.01$. {\bf (B)} A 3-dimensional perspective of (A) with the same parameters, for better understanding. {\bf (C)} Three simulations are done, with same parameters and $\gamma_i=0.05$, for three distinct prevalence thresholds $P=0.01,\ 0.05,\ 0.1$ .}
\label{fig:EFS_Information_R0_MultipleShot}
\end{figure}

\begin{figure}[H]
\centering
\includegraphics[width=\textwidth]{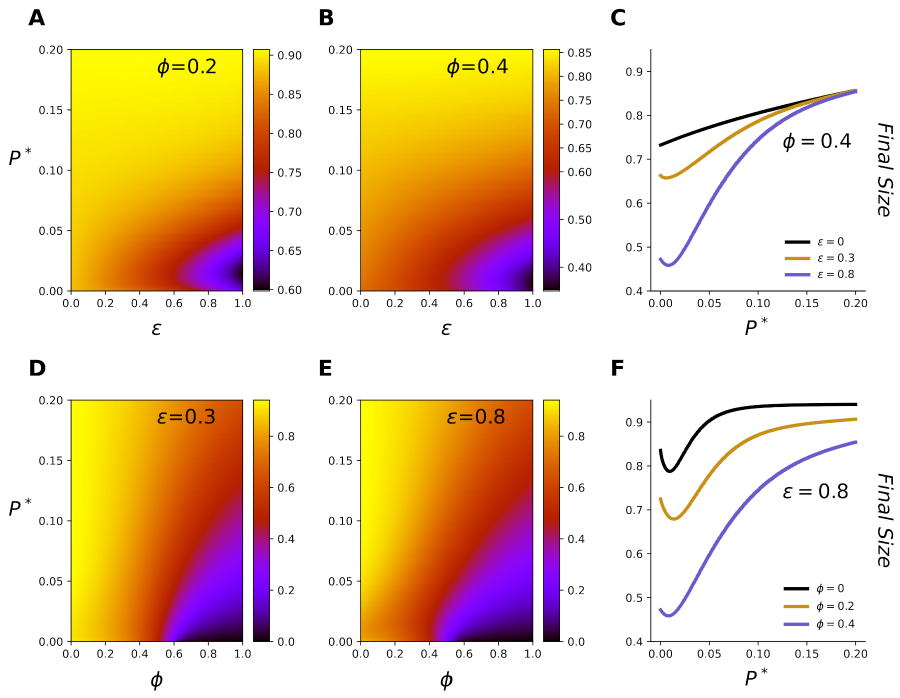}
\caption{{\bf Change in epidemic final size by varying the prevalence threshold and behavioral parameters in a multiple shot framework.} {\bf (A and B)} Different scenarios by modifying the parameter $\phi=0.2,\ 0.4$ respectively for (A) and (B) to get the epidemic final size, where $\epsilon$ is varied between $[0,1]$ and $P$ between $[0,1]$, while keeping the rest of the parameters fixed, $\rho=1/20$ (multiple shot), $\beta =0.6$, $\gamma =0.2$, $\beta_i=1.5$ and $\gamma_i=0.1$. {\bf (D and E)} Images representing the same as in (A and B), with the difference being that the epidemic final size is obtained by varying $\phi$ instead of $\epsilon$, and considering two cases for $\epsilon=0.3,\ 0.8$ respectively for (D and E). {\bf (C and F)} Vertical trajectories obtained from (B and E) for (C and F) respectively.}
\label{fig:BehaviorParameters_grids_MultipleShot}
\end{figure}

% \begin{figure}[H]
% \centering
% \includegraphics[width=\textwidth]{Comparison_Behavioral_Params_MultipleShot.pdf}
% \caption{{\bf Comparison between behavioral response strength and enforced behavioral response in a multiple shot configuration.} {\bf (A and B)} Representation of the epidemic final size based on the comparison between behavioral parameters. A contour plot is used, with parameters such that $R_0^d = 6$ for (A) and $R_0^d = 2.22$ for (B), respectively. Multiple shot configuration is used.}
% \label{fig:Comparison_Behavioral_Params_MultipleShot}
% \end{figure}

%%%%%%%%%%%%%%%% SUPPLEMENTARY TABLES %%%%%%%%%%%%%%%

%%%%%%%%%%% CAPTIONS FOR OTHER SUPPLEMENTARY FILES %%%%%%%%%%

\clearpage % Clear all remaining figures and tables then start a new page

%%%%%%%%%%%%%%%% SUPPLEMENTARY REFERENCES %%%%%%%%%%%%%%%

% Do NOT include a reference list in the supplement.
% All references must be in a single list at the end of the main text.
% The copyeditors will ensure that the correct reference list appears with each version of the paper
% (print, HTML, PDF, mobile app, metadata for bibliographic databases etc.)

\end{document}